\documentclass[12pt]{article}
\usepackage{amsmath,amsthm, amsfonts, amssymb, latexsym,verbatim,bbm,longtable}
\input xy
\xyoption{all}
\usepackage{hyperref,multirow}
\usepackage{tikz, ifthen}
\usepackage[shortlabels]{enumitem}
\usepackage{young}
\usepackage{mathtools}
\usepackage{color,soul}
\usepackage{dynkin-diagrams}
\usepackage{standalone}
\usepackage{multicol} 
\usepackage{caption}
\usepackage{subcaption}
\usepackage{float}
\usepackage{dynkin-diagrams}
\usepackage{stmaryrd}
\usepackage{graphicx}

\DeclareMathOperator{\gap}{gap}
\DeclareMathOperator{\maxx}{max}
\DeclareMathOperator{\minn}{min}
\DeclareMathOperator{\uc}{u}
\DeclareMathOperator{\dc}{d}
\DeclareMathOperator{\edges}{edges}
\DeclareMathOperator{\maxchains}{c}

\newcommand{\R}{\mathbb{R}}

\newcommand{\OO}{\mathcal{O}}
\newcommand{\CC}{\mathcal{C}}

\newcommand{\newword}{\textbf}

\theoremstyle{plain}
\newtheorem{theorem}{Theorem}
\newtheorem{lemma}[theorem]{Lemma}

\theoremstyle{definition}
\newtheorem{definition}[theorem]{Definition}
\newtheorem{example}[theorem]{Example}

\newtheorem{question}[theorem]{Question}
\newtheorem{remark}[theorem]{Remark}

\title{Gap between the number of facets of the two poset polytopes}

\author{Binaya Bhandari \thanks{Email: cqy14@txstate.edu} \and Debra Cunningham \thanks{Email: db01609@txstate.edu} \and Grace Morrell \thanks{Email: efj19@txstate.edu} \and SuHo Oh \thanks{Email: s\_o79@txstate.edu} \and Paxton Smith \thanks{Email: jps184@txstate.edu}}

\begin{document}
\maketitle

\begin{abstract}
    We study the difference between the number of facets of the order polytope and the chain polytope of a poset. Hibi and Li classified posets where the gap is exactly zero. We describe the bounds on this gap using the new notion of crossing numbers, and then use this result to classify the posets where the gap is exactly one.
\end{abstract}

\section{Introduction}

Let $P$ be a finite poset with $n$ elements. Stanley in \cite{stanley1986two} defined the following two polytopes:

\begin{definition}[\cite{stanley1986two}]
 The \newword{order polytope} of a poset $P$ is given as
  $$\OO(P) = \{x \in \R^n_{\geq 0} \mid x_i \leq 1 \text{ and } x_i \leq x_j \text{ if } i \leq j \text{ in $P$}\}.$$ 
 The \newword{chain polytope} of a poset $P$ is given as
  $$\CC(P) = \{x \in \R^n_{\geq 0} \mid x_i \leq 1 \text{ and } x_{p_1} + \dots + x_{p_s} \leq 1 \text{ for all } p_1 \leq \dots \leq p_s \text{ in $P$} \}.$$
  \end{definition}
  

Both polytopes have full dimension $n$, and are known to have the same Ehrart polynomial \cite{stanley1986two}. So, it is natural to compare these two similar polytopes.

Recall that the \newword{$f$-vector} $(f_0,\ldots,f_{n-1})$ of an $n$-dimensional polytope counts the number of faces per dimension. The number of vertices is known to be the same between the two polytopes ($f_0(\OO(P) = f_0(\CC(P))$) \cite{stanley1986two}. The number of facets of the order polytope is always less than or equal to the number of facets of the chain polytope ($f_{n-1}(\OO(P)) \leq f_{n-1}(\CC(P))$) \cite{HIBILI16}, and a similar result holds for the number of edges as well ($f_{1}(\OO(P)) \leq f_{1}(\CC(P))$) \cite{HIBI17}. In \cite{HIBILI16}, it was conjectured that this inequality extends to any other $f_i$. This conjecture was partially resolved for a specific class of ranked posets called \newword{maximal ranked posets} in \cite{ahmad2023order}. In \cite{freijhollanti2024fvector}, it was shown that the posets that come from disjoint unions and ordinal sums of posets that already satisfy the inequality again satisfy the inequality.

In this paper, we focus on comparing the number of facets of the two polytopes for general posets. Hibi and Li in \cite{HIBILI16}, in addition to showing $f_{n-1}(\OO(P)) \leq f_{n-1}(\CC(P))$, also showed the following theorem, which describes exactly when the two polytopes have the same number of facets:

\begin{theorem}[\cite{HIBILI16}]
\label{thm:hibiX}
We have $f_{n-1}(\OO(P)) = f_{n-1}(\CC(P))$ if and only if the poset $P$ does not contain five pairwise distinct elements $a,b,c,d,e$ such that $a,b \prec c \prec d,e$.
\end{theorem}

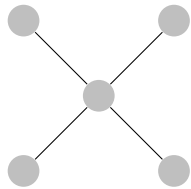
\begin{figure}[h]

 \begin{center}
 
  \vspace{3.7mm}
    \begin{tikzpicture}[shorten >=1pt,->]
  \tikzstyle{vertex}=[circle,fill=black!25,minimum size=12pt,inner sep=2pt]
  \node[vertex] (G_1) at (-1,-1) {};
  \node[vertex] (G_2) at (0,0)   {};
  \node[vertex] (G_3) at (1,-1)  {};
  \node[vertex] (G_4) at (-1,1)  {};
  \node[vertex] (G_5) at (1,1)   {};
  \draw (G_1) -- (G_2) -- (G_3) -- cycle;
  \draw (G_4) -- (G_2) -- (G_5) -- cycle;
\end{tikzpicture}
    \captionsetup{width=1.0\linewidth}
  \captionof{figure}{The $X$-shaped poset.}
  \label{fig:x-shape}

 \end{center}

\end{figure}

The poset described in Theorem~\ref{thm:hibiX} is drawn in Figure~\ref{fig:x-shape}. We call such a poset an \newword{X-shaped} poset and call $c$ its \newword{center}. If a poset does not contain an X-shaped subposet, we say that it is \newword{X-avoiding}.

In this paper, we are going to focus on extending this result. Our main idea comes from defining a new statistic on the vertices of a poset, called \newword{crossing numbers}. This allows us to arrive at the bounds of $f_{n-1}(\CC(P)) - f_{n-1}(\OO(P))$ (which we call \newword{gap} of $P$). In Section $2$, we introduce our main result, the bound on the gap. 
In Section $3$, using this bound, we classify the posets where the gap is exactly one.

\section{The bound on the gap}

We start by defining the gap of a poset $P$:

\begin{definition}
Given a poset $P$ on $n$ elements, we will denote $f_{n-1}(\CC(P)) - f_{n-1}(\OO(P))$, the difference between the number of facets of the chain polytope and the order polytope, as \newword{gap} of $P$ and write it as $\gap(P)$. 
\end{definition}

For any poset $P$, we have $\gap(P) \geq 0$ from Corollary 1.2 of \cite{HIBILI16}. We now review how to obtain $f_{n-1}(\CC(P))$ and $f_{n-1}(\OO(P))$ from the poset $P$. In this paper, we use $\edges(P)$ to denote the number of edges in the Hasse diagram of $P$, use $\maxchains(P)$ to denote the number of maximal chains of $P$, use $n$ to represent the number of elements of $P$, use $\maxx(P)$ to denote the number of maximal elements of $P$ and use $\minn(P)$ to denote the number of minimal elements of $P$. 

\begin{lemma}[Lemma 1.1 of \cite{HIBILI16}]
\label{lem:facets}
The number of facets of the order polytope of $P$ is given by $\maxx(P) + \minn(P) + \edges(P)$ and the number of facets of the chain polytope is given by $\maxchains(P) + |P|$. 
\end{lemma}

Below are examples of calculating the \newword{gap} using Lemma~\ref{lem:facets}.

\begin{example}
    Consider the poset in Figure 2. Notice that this poset is X-avoiding. 
    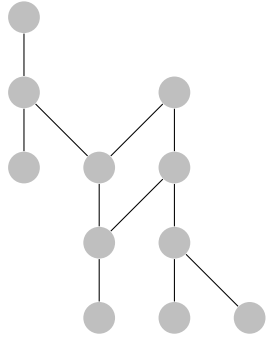
\begin{figure}[h]

 \begin{center}
 
  \vspace{3.7mm}
    \begin{tikzpicture}[shorten >=1pt,->]
  \tikzstyle{vertex}=[circle,fill=black!25,minimum size=12pt,inner sep=2pt]
  \node[vertex] (J_1) at (-1,2) {};
  \node[vertex] (J_2) at (-1,1)   {};
  \node[vertex] (J_3) at (-1,0)  {};
  \node[vertex] (J_4) at (0,0)  {};
  \node[vertex] (J_5) at (0,-1)   {};
  \node[vertex] (J_6) at (0,-2)   {};
  \node[vertex] (J_7) at (1,1)   {};
  \node[vertex] (J_8) at (1,0)   {};
  \node[vertex] (J_9) at (1,-1)   {};
  \node[vertex] (J_10) at (1,-2)   {};
  \node[vertex] (J_11) at (2,-2)   {};

  \draw (J_1) -- (J_2) -- (J_3) -- cycle;
  \draw (J_2) -- (J_4) -- (J_5) -- (J_6) -- cycle;
  \draw (J_4) -- (J_7) -- (J_8) --  (J_9) --  (J_10) -- cycle;  
  \draw (J_9) -- (J_11) -- cycle;
  \draw (J_5) -- (J_8) -- cycle;
\end{tikzpicture}
    \captionsetup{width=1.0\linewidth}
  \captionof{figure}{X-avoiding poset}
  \label{fig:non-x-shape ex 2}

 \end{center}
\end{figure}
\begin{align*}
\gap(P) & = c(P) + |P| - (\maxx(P)+ \minn(P)+ \edges(P))\\
&= 6 + 11 - (2+4+11)\\
&=0
\end{align*}
\end{example}

\begin{example}
    Consider the poset in Figure 3. Notice that this poset contains an X-shaped subposet.  
    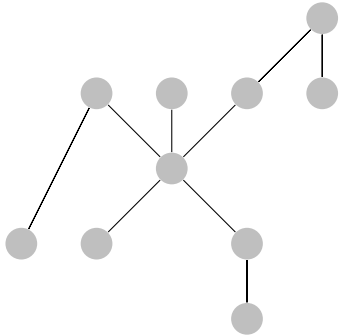
\begin{figure}[h]

 \begin{center}
 
  \vspace{3.7mm}
    \begin{tikzpicture}[shorten >=1pt,->]
  \tikzstyle{vertex}=[circle,fill=black!25,minimum size=12pt,inner sep=2pt]
  \node[vertex] (G_1) at (-1,-1) {};
  \node[vertex] (G_2) at (0,0)   {};
  \node[vertex] (G_3) at (1,-1)  {};
  \node[vertex] (G_4) at (-1,1)  {};
  \node[vertex] (G_5) at (1,1)   {};
  \node[vertex] (G_6) at (2,1)   {};
  \node[vertex] (G_7) at (1,-2)   {};
  \node[vertex] (G_8) at (2,2)   {};
  \node[vertex] (G_9) at (-2,-1)   {};
  \node[vertex] (G_10) at (0,1)   {};
  \draw (G_1) -- (G_2) -- (G_3) -- cycle;
  \draw (G_4) -- (G_2) -- (G_5) -- cycle;
  \draw (G_5) -- (G_8) -- (G_5) -- cycle;
  \draw (G_8) -- (G_6) -- (G_8) -- cycle;
  \draw (G_3) -- (G_7) -- (G_3) -- cycle;
  \draw (G_9) -- (G_4) -- (G_9) -- cycle;
  \draw (G_2) -- (G_10) -- cycle;
\end{tikzpicture}
    \captionsetup{width=1.0\linewidth}
  \captionof{figure}{Poset with an X-shaped subposet}
  \label{fig:non-x-shape ex 3}

 \end{center}
\end{figure}
\begin{align*}
\text{gap(P)} & = c(P) + |P| - (\maxx(P)+ \minn(P)+ \edges(P))\\
&= 8 + 10 - (3+4+9)\\
&=2
\end{align*}
\end{example}


 



In the proof of Theorem~\ref{thm:hibiX} in \cite{HIBILI16}, the main idea was to look at the difference between $gap(P)$ and $gap(P \setminus \alpha)$ where $\alpha$ is a minimal element. Let $\beta_1,\ldots,\beta_s,\gamma_1,\ldots,\gamma_t$ denote the elements of $P$ that cover $\alpha$, such that $\beta_i$ covers at least two elements of $P$ and each $\gamma_j$ covers no element of $P$ except for $\alpha$. That is, $\gamma_j$'s are the elements that become minimal in $P \setminus \alpha$ and $\beta_i$'s are the elements that are still not minimal in $P \setminus \alpha$. For example, look at Figure~\ref{fig:image1}.

\begin{definition}
    Let $v$ be any element of the poset $P$. We call any saturated chain having some maximal element of $P$ as its maximal element, and having $v$ as its minimal element to be a \newword{maximal chain above $v$}. Similarly, we call any saturated chain having some minimal element of $P$ as its minimal element, and having $v$ as its maximal element to be a \newword{maximal chain below $v$}. We use $\uc(v)$ to denote the number of maximal chains above $v$ and $\dc(v)$ to denote the number of maximal chains below $v$.
\end{definition}

The following lemma is given implicitly in the proof of Corollary 1.2 and Theorem 1.3 of \cite{HIBILI16}, but we put the full proof here as well for the convenience of the reader.

\begin{lemma}[\cite{HIBILI16}]
\label{lem:gapdiff}
    Let $\alpha$ be any minimal and nonmaximal element of $P$. Then we have the following:
    $$\gap(P) - \gap(P \setminus \alpha) = \sum_{i=1}^s \uc(\beta_i)- s.$$
\end{lemma}
\begin{proof}
Comparing $P$ to $P \setminus \alpha$, $P$ has $1$ more node, $s+t$ more edges, and $(t-1)$ less minimal elements. We analyze the difference between $\maxchains(P)$ and $\maxchains(P \setminus \alpha)$.

Any maximal chain of $P$ that does not have $\alpha$ as its minimal element is again a maximal chain in $P \setminus \alpha$. We will create an injective map $\phi$ from the set of maximal chains of $P \setminus \alpha$ to the set of maximal chains of $P$. If a maximal chain $C$ of $P \setminus \alpha$ does not have any of $\gamma_1,\ldots,\gamma_t$ as its minimal element, let $\phi(C)$ again be $C$. If a maximal chain $C$ of $P \setminus \alpha$ has one of $\gamma_1,\ldots,\gamma_t$ as its minimal element, let $\phi(C)$ be the maximal chain obtained from $C$ by attaching $\alpha$. Then the maximal chains of $P$ that are not in the image of $\phi$ are exactly those that have $\alpha$ as its minimal element and go through one of $\beta_1,\ldots,\beta_s$. This set is counted by $\sum_{i=1}^s \uc(\beta_i)$.
\begin{multline*}
\gap(P) - \gap(P \setminus \alpha) = 
\maxchains(P) - \maxchains(P \setminus \alpha) + |P| - |P \setminus \alpha| \\ - (\maxx(P) - \maxx(P \setminus \alpha) + \minn(P) - \minn(P \setminus \alpha) + 
\edges(P) - \edges(P \setminus \alpha)) \\
= \sum_{i=1}^s \uc(\beta_i) + 1 - (0 - (t-1) + t+s)
= \sum_{i=1}^s \uc(\beta_i) - s
\end{multline*}

\end{proof}

Notice that the right-hand side of the above equation can also be written as $\sum_{i=1}^s(\uc(\beta_i)-1)$. Take a look at Example~\ref{ex:calcgap} on how to use Lemma~\ref{lem:gapdiff} to calculate the difference between the gaps.

\begin{example}
\label{ex:calcgap}
    Consider the poset in Figure 4. The difference between $\gap(P)$ and $\gap(P \setminus \alpha)$ is given by $\sum_{i=1}^2 \uc(\beta_i)-2$. Since $\uc(\beta_1) = \uc(\beta_2) = 2$, we have that the difference between the gaps is $2$.
    \begin{figure}[h]
    \centering
    \begin{subfigure}{0.33\linewidth}
        \centering
        \begin{tikzpicture}[> = stealth,  shorten > = 1pt,   auto,   node distance = 1.5cm,
every edge/.style = {draw=black,thick},
 vrtx/.style args = {#1/#2}{%
      circle, draw, thick, fill=black!25,
      minimum size=5mm, label=#1:#2}
                ]
\node (A) [vrtx=left/{}]     at ( 0, 0) {};
\node (B) [vrtx=left/{}]     at (-1, 1) {};
\node (C) [vrtx=right/{}]    at ( 1, 1) {};
\node (D) [vrtx=left/$\beta_1$]    at ( 1.1,-1) {};
\node (E) [vrtx=left/$\gamma_1$]    at (-1.5,-1) {};
\node (F) [vrtx=left/$\alpha$]    at (1,-2) {};
\node (G) [vrtx=left/$\gamma_2$]    at (-0.1,-1) {};
\node (H) [vrtx=left/$\beta_2$]    at (2.2,-1) {};
\node (I) [vrtx=left/{}]    at (2.5,-2) {};
\node (J) [vrtx=left/{}]    at (4,-2) {};
 \path   (A) edge (B)
         (A) edge (C)
       (E) edge (A)
         (G) edge (A)
         (D) edge (A)
         (H) edge (A)
         (F) edge (E)
         (F) edge (G)
         (F) edge (H)
         (D) edge (I)
         (H) edge (J)
         (F) edge (D);
\end{tikzpicture}
        \caption{A poset $P$ with a minimal element $\alpha$}
        \label{fig:image1}
    \end{subfigure}
    \hskip 3cm
    \begin{subfigure}{0.33\linewidth}
        \centering
        \begin{tikzpicture}[> = stealth,  shorten > = 1pt,   auto,   node distance = 1.5cm,
every edge/.style = {draw=black,thick},
 vrtx/.style args = {#1/#2}{%
      circle, draw, thick, fill=black!25,
      minimum size=5mm, label=#1:#2}
                    ]
\node (A) [vrtx=left/{}]     at ( 0, 0) {};
\node (B) [vrtx=left/{}]     at (-1, 1) {};
\node (C) [vrtx=right/{}]    at ( 1, 1) {};
\node (D) [vrtx=left/$\beta_1$]    at ( 1.1,-1) {};
\node (E) [vrtx=left/$\gamma_1$]    at (-1.5,-1) {};
\node (G) [vrtx=left/$\gamma_2$]    at (-0.1,-1) {};
\node (H) [vrtx=left/$\beta_2$]    at (2.2,-1) {};
\node (I) [vrtx=left/{}]    at (2.5,-2) {};
\node (J) [vrtx=left/{}]    at (4,-2) {};
 \path   (A) edge (B)
         (A) edge (C)
       (E) edge (A)
         (G) edge (A)
         (D) edge (A)
         (H) edge (A)
         (D) edge (I)
         (H) edge (J);
\end{tikzpicture}
         \caption{The poset $P \setminus \alpha$}
        \label{fig:image2.pdf}
    \end{subfigure}
    \caption{How the labeling of $\alpha,\beta_1,\ldots,\beta_s,\gamma_1,\ldots,\gamma_t$ works}
\end{figure}
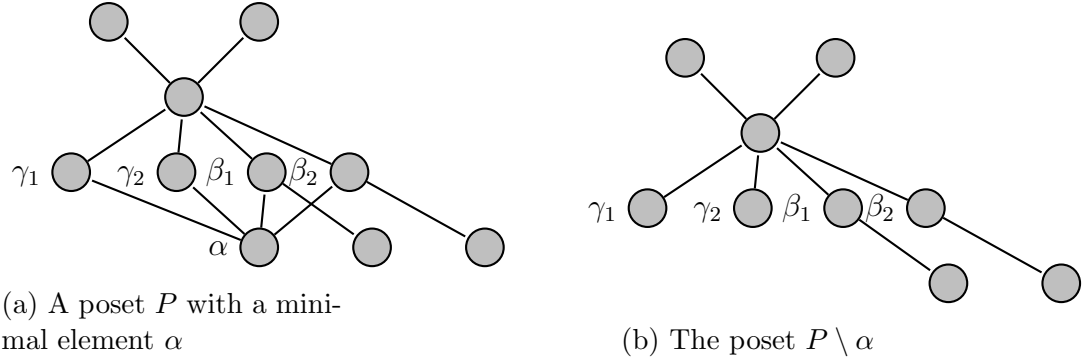


 

\end{example}

Now we state and prove our main result. Let \newword{crossing number} of $v \in P$ be $(\uc(v)-1)(\dc(v)-1)$. Notice that $\uc(v)\dc(v)$ is simply the number of maximal chains of $P$ that pass through $v$. If for any maximal antichain $A$ we take $\sum_{v \in A} \uc(v)\dc(v)$, we get the number of maximal chains of $P$.

\begin{theorem}
\label{thm:main}
Let $P$ be a poset and $A$ be any maximal antichain of $P$. Then we have the following:

$$ \sum_{v \in A} (\uc(v) - 1)(\dc(v)-1) \leq \gap(P) \leq \sum_{v \in A} (\uc(v)\dc(v)-1).$$
i.e. the gap is bounded below by the sum of crossing numbers for any maximal antichain $A$ of $P$.
\end{theorem}

\begin{proof}

We are going to use induction on the number of nodes to prove both the lower bound and the upper bound on the gap. We start with the base case: If all the minimal elements of $P$ are maximal elements as well, then $P$ is simply an antichain itself. In this case, the gap is $0$ due to Theorem~\ref{thm:hibiX}. All vertices have a crossing number of $0$ so the left-hand side is zero, and all vertices have exactly one maximal chain passing through them, so the right-hand side is zero as well.

Now assume for the sake of induction that we have the claim for posets having up to $n-1$ vertices, and consider the poset $P$ that has $n$ vertices. We track the difference in the lower bound, the gap, and the upper bound between $P$ and $P \setminus \alpha$ when $\alpha$ is a minimal and nonmaximal element of $P$. Recall that from Lemma~\ref{lem:gapdiff}, we have $\gap(P) - \gap(P \setminus \alpha) = \sum_{i=1}^s \uc(\beta_i) - s$.

First, look at the upper bound, when $\alpha \in A$. Let $A'$ be obtained from $A$ by removing $\alpha$ and adding all elements among $\beta_1,\ldots,\beta_s,\gamma_1,\ldots,\gamma_t$ that are incomparable with elements of $A \setminus \alpha$ (there has to be at least one, since otherwise $A$ is not a maximal antichain of $P$). When $\alpha \not \in A$, simply set $A'$ as $A$.

Then in either case, $A'$ is a maximal antichain in both $P$ and $P \setminus \alpha$ such that $-|A'| \leq - |A|$. Notice that the upper bound is simply the number of maximal chains of $P$ minus the cardinality of $A$. Since $-|A'| \leq -|A|$, it is enough to prove the claim for $A'$ instead. The difference between $P$ and $P \setminus \alpha$ of the value $\maxchains(P) - |A'|$ is exactly the number of maximal chains that go through one of $\beta_i$ and end at $\alpha$, which is given by $\sum_{i=1}^s \uc(\beta_i)$. This is clearly greater than or equal to $\sum_{i=1}^s \uc(\beta_i) - s$. Hence, by induction, we have proved the upper bound on the gap.


We now look at the lower bound. Recall that from the induction hypothesis, we have $\sum_{v \in A}(\uc(v)-1)(\dc(v)-1) \leq \gap(P \setminus \alpha)$ for any maximal antichain $A$ of $P \setminus \alpha$. Now in order to show that the lower bound holds for $P$, since $\alpha$ has the crossing number of zero, it is enough to prove that the same inequality holds for any maximal antichain of $P$ that does not contain $\alpha$: that is, for any maximal antichain of $P \setminus \alpha$.

So, let $A$ be a maximal antichain of $P \setminus \alpha$. We let $A' = \{a_1,\ldots,a_q\}$ be the elements of $A$ that is greater than or equal to some $\beta_i$. Relabel the elements $\beta_1,\ldots,\beta_s$ so that for all $\beta_i$, there exists some $s' \leq s$ such that $\beta_i$ for $i \leq s'$ is less or equal to some element of $A'$ and $\beta_i$ for $i > s'$ is not.

We use $y_{i,j}$ to denote the number of saturated chains that have $\beta_i$ ($1 \leq i \leq s'$) as its minimal element and $a_j$ ($1 \leq j \leq q$) as its maximal element. Then the difference of the sum $\sum_{v \in A}(\uc(v)-1)(\dc(v)-1)$ as we go from $P \setminus \alpha$ to $P$ is given by $\sum_{j=1}^q (\uc(a_j)-1) \sum_{i=1}^{s'} y_{i,j}$, since the number of maximal chains above $a_j$ does not change and the number of maximal chains below $a_j$ increases exactly by $\sum_{i=1}^{s'} y_{i,j}$. The above difference can be rewritten as $\sum_{i=1}^{s'} \sum_{j=1}^q \uc(a_j) y_{i,j} - \sum_{i=1}^{s'}\sum_{j=1}^q y_{i,j}$. 

From the definition of $A'$ and $s'$, we have $\sum_{j=1}^{q} y_{i,j} \geq 1$ for each $i \leq s'$, which gives us $\sum_{i=1}^{s'}\sum_{j=1}^q y_{i,j} \geq s'$. From the definition of $y_{i,j}$'s we have $\sum_{j=1}^q \uc(a_j)y_{i,j} \leq \uc(\beta_i)$. So we see that the previous difference is bounded above by $\sum_{i=1}^{s'} (\uc(\beta_i)-1)$. Since $\uc(\beta_i) \geq 1$ for each $1 \leq i \leq s$, we see that this is again bounded above by $\sum_{i=1}^s (\uc(\beta_i)-1)$, which is exactly the difference in the gap as we go from $P \setminus \alpha$ to $P$ thanks to Lemma~\ref{lem:gapdiff}. By induction, this finishes the proof of the lower bound.
   
\end{proof}

\begin{remark}
The bounds of Theorem~\ref{thm:main} are tight in the sense that there exist some posets where both bounds collapse to the same number and become an equality overall. In particular, consider a chain of any length. The gap is obviously $0$ from Theorem~\ref{thm:hibiX}. Pick any vertex $v$, which is by itself a maximal antichain, and the value $\uc(v)\dc(v)-1$ is zero, since there is only one unique maximal chain of this poset. The crossing number for any vertex in a chain is zero. So both the upper bound and the lower bound give zero here.
\end{remark}



\section{Classifying posets that have gap being $1$}
In this section, using our main result of Theorem~\ref{thm:main}, we will classify the posets that have the gap to be exactly $1$. We first start by discussing how the X-avoiding posets can be described using the crossing numbers:

\begin{theorem}
\label{thm:zerogap}
Let $P$ be a finite poset. Then the following are equivalent:
\begin{enumerate}
    \item $P$ is $X$-avoiding.
    \item gap of $P$ is $0$.
    \item for all $v \in P$, the crossing number of $v$ is zero.
\end{enumerate}

\end{theorem}
\begin{proof}
Recall that from Theorem~\ref{thm:hibiX}, we see that $(1)$ and $(2)$ are equivalent. Now, from Theorem~\ref{thm:main}, the gap is bounded below by the sum of crossing numbers for any antichain of $P$. So, this means that if the gap of $P$ is zero, then all vertices must also have the crossing number of zero. So we have shown that $(2) \rightarrow (3)$.

Finally, we are going to show that $(3) \rightarrow (1)$. When for all $v \in P$ the crossing number is zero, assume that for the sake of contradiction we have some X-shaped subposet inside our $P$ with $v$ as its center. Then, since $\uc(v) \geq 2$ and $\dc(v) \geq 2$, the crossing number of $v$ is at least one, which gives us a contradiction.
\end{proof}

Now we classify the posets where the gap is one. To do so, we first define the following.

\begin{definition}
\label{def:xorchid}
We say that a finite poset $P$ is an \newword{X-orchid} if satisfies the following conditions. It has a saturated chain $c_1 < \dots < c_k$ where $c_k$ is covered above by exactly two elements $u_1,u_2$ and $c_1$ is covered below by exactly two elements $d_1,d_2$. Moreover, any subposet that is X-shaped must have some $c_i$ as its center, and by replacing $c_i$ with any other $c_j$, we get another X-shaped subposet inside $P$ with $c_j$ as its center.
\end{definition}

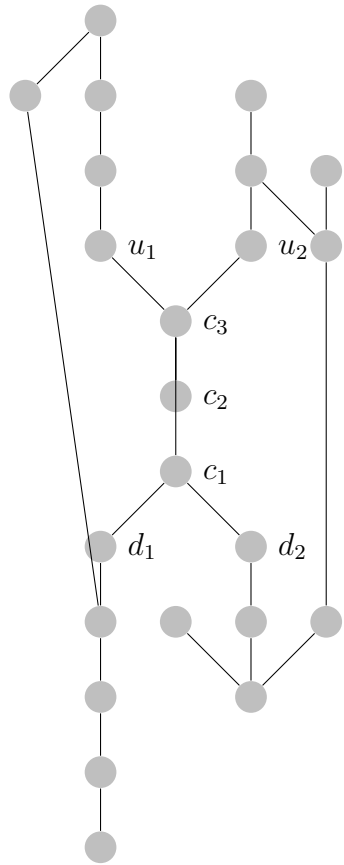
\begin{figure}[h]

 \begin{center}
 
  \vspace{3.7mm}
    \begin{tikzpicture}[shorten >=1pt,->]
  \tikzstyle{vertex}=[circle,fill=black!25,minimum size=12pt,inner sep=2pt]
  \node[vertex, label=right:{\textcolor{black}{$c_2$}}] (P_1) at (0,0) {};
  \node[vertex, label=right:{\textcolor{black}{$c_3$}}] (P_2) at (0,1) {};
  \node[vertex, label=right:{\textcolor{black}{$c_1$}}] (P_5) at (0,-1) {};
  \node[vertex, label=right:{\textcolor{black}{$u_2$}}] (P_3) at (1,2) {};
  \node[vertex] (P_3a) at (1,3) {};
  \node[vertex] (P_3b) at (2,2) {};
  \node[vertex] (P_3c) at (1,4) {};
  \node[vertex] (P_3d) at (2,3) {};
  \node[vertex, label=right:{\textcolor{black}{$u_1$}}] (P_4) at (-1,2) {};
  \node[vertex] (P_4a) at (-1,3) {};
  \node[vertex] (P_4b) at (-1,4) {};
  \node[vertex] (P_4c) at (-1,5) {};
  \node[vertex] (P_4d) at (-2,4) {};
  \node[vertex, label=right:{\textcolor{black}{$d_2$}}] (P_6) at (1,-2) {};
  \node[vertex] (P_6a) at (1,-3) {};
  \node[vertex] (P_6b) at (1,-4) {};
  \node[vertex] (P_6c) at (0,-3) {};
  \node[vertex] (P_6d) at (2,-3) {};
  \node[vertex, label=right:{\textcolor{black}{$d_1$}}] (P_7) at (-1,-2) {};
  \node[vertex] (P_7a) at (-1,-3) {};
  \node[vertex] (P_7b) at (-1,-4) {};
  \node[vertex] (P_7c) at (-1,-5) {};
  \node[vertex] (P_7d) at (-1,-6) {};

  \draw (P_1) -- (P_2) -- (P_5) -- cycle;
  \draw (P_3b) -- (P_6d) -- cycle;
  \draw (P_4d) -- (P_7a) -- cycle;
  \draw (P_3) -- (P_2) -- (P_4) -- cycle;
  \draw (P_6) -- (P_5) -- (P_7) -- cycle;  
  \draw (P_6) -- (P_6a) -- (P_6b) -- cycle;
  \draw (P_6c) -- (P_6b) -- (P_6d) -- cycle;
  \draw (P_3) -- (P_3a) -- (P_3c) -- cycle;
  \draw (P_3a) -- (P_3b) -- (P_3d) -- cycle;
  \draw (P_4) -- (P_4a) -- (P_4b) -- (P_4c) -- (P_4d) -- cycle;
  \draw (P_7) -- (P_7a) -- (P_7b) -- (P_7c) -- (P_7d) -- cycle;
  
\end{tikzpicture}
    \captionsetup{width=1.0\linewidth}
  \captionof{figure}{Example of an X-Orchid}
  \label{fig:X-Orchid}

 \end{center}

\end{figure}

An example of an X-orchid is shown in Figure~\ref{fig:X-Orchid}. 
Notice that $c_3$ is covered above by exactly two elements, and $c_1$ covers exactly two elements below. We call the saturated chain $c_1 < \dots < c_k$ together with $u_1,u_2,d_1,d_2$ the \newword{stalk} of the X-orchid.

\begin{theorem}
\label{thm:gapone}
Let $P$ be a finite poset. Then the following are equivalent:
\begin{enumerate}
    \item $P$ is an X-orchid.
    \item gap of $P$ is $1$.
    \item $P$ has a saturated chain consisting of elements that have a crossing number of one. All other elements have zero crossing number.
\end{enumerate}
\end{theorem}
\begin{proof}
We may assume that $P$ does not have any element that is minimal and maximal at the same time, since such elements have a zero crossing number and also do not affect the gap.

We begin by showing $(1) \rightarrow (2)$. We show this by induction on the number of elements. We start with the base case: Let $P$ be an X-orchid that consists only of its stalk: a saturated chain $c_1 < \dots < c_k$ with $u_1,u_2$ covering $c_k$, and $c_1$ covering $d_1,d_2$. Then we have $\gap(P) = 4 + (k+4) - (2 + 2 + k+3) = 1$. 

Next, we consider the case where $d_1,d_2$ of the stalk are the only minimal elements of $P$. Given any X-shaped subposet of $P$, we can replace its center with $c_1$ to get another X-shaped subposet due to the definition of an X-orchid. This means that $P \setminus d_1$ is X-avoiding. So, we have $\gap(P \setminus d_1) = 0$ thanks to Theorem~\ref{thm:zerogap}. Setting $\alpha = d_1$, using the notation from \cite{HIBILI16}, let $\beta_1,\ldots,\beta_s,\gamma_1,\ldots,\gamma_t$ denote the elements that cover $\alpha$. Since the only minimal elements of $P$ are $d_1$ and $d_2$, we have $s=1$ and $\beta_1 = c_1$. Using Lemma~\ref{lem:gapdiff}, we get $\gap(P) - \gap(P \setminus d_1) = \uc(c_1)-1$. If for some $v > c_k$ we have at least two elements that cover $v$, then we get an X-shaped subposet with $v$ as the center, so this cannot happen. This implies $\uc(c_1) = 2$, giving us that $\gap(P) = 1$.

Assume for the sake of induction that we have the claim that $P$ has strictly fewer vertices than $n$. Let $P$ be an X-orchid with $n$ vertices that contains some minimal element $\alpha \not \in \{d_1,d_2\}$.  Again using the notation from \cite{HIBILI16}, let $\beta_1,\ldots,\beta_s,\gamma_1,\ldots,\gamma_t$ denote the elements of $P$ that cover $\alpha$, so that $\beta_i$ covers at least two elements. Then from the definition of an X-orchid, we have $\uc(\beta_i) = 1$ since otherwise we can find some X-shaped subposet having $\beta_i$ as its center. From Lemma~\ref{lem:gapdiff}, we have $\gap(P) = \gap(P \setminus \alpha)$. Hence, by induction we have shown that $\gap(P) = 1$.


Now we focus on showing $(2) \rightarrow (3)$. From Theorem~\ref{thm:main}, the crossing number of each $v \in P$ must be either zero or one. Moreover, there cannot be pairwise incomparable elements that both have a nonzero crossing number. Also, from Theorem~\ref{thm:zerogap}, there has to exist some $v \in P$ with the crossing numer being one.

Now, if we have $x < y < z$ in $P$ where $x$ and $z$ both have a crossing number of one, then $y$ cannot have a crossing number of zero. This is because $\uc(y) \geq \uc(z)$ and $\dc(y) \geq \dc(x)$ and we know that $\uc(x), \dc(z) > 1$ from the fact that $x,z$ has a non-zero crossing number. So, the set of $v \in P$ with crossing number equal to one has to form a saturated chain in $P$ and all other vertices have crossing number of zero. 

Finally, we focus on showing $(3) \rightarrow (1)$. Let $c_1 < \dots < c_k$ be the saturated chain of elements all having crossing number one. Let $u_1,\ldots,u_r$ be the elements that cover $c_k$ and let $d_1,\ldots,d_s$ be the elements that are covered by $c_1$. Since the crossing number of $c_i$ is one, we must have $r=s=2$. Moreover, each $c_i$ for $i < k$ is covered above by exactly one element, and each $c_i$ for $i > 1$ covers exactly one element below. So, if there is an X-shaped subposet that has $c_i$ as its center, we can replace $c_i$ with any other $c_j$ to get an X-shaped subposet within $P$.

Now assume for the sake of contradiction that $P$ contains some X-shaped subposet, where the center is not one of $c_i$'s. The center of this X-shaped subposet must have a crossing number of at least one, so we get a contradiction. 
\end{proof}

So we have classified the posets that have the gap as exactly one. From the above theorem, it seems very natural to ask the following questions:

\begin{question}
    Can one classify the posets having gap of two using similar methods?
\end{question}

\begin{question}
    Is it possible to come up with a better upper bound for Theorem~\ref{thm:main} that gives the value of $0$ for X-avoiding posets and $1$ for X-orchids to make the proof of Theorem~\ref{thm:zerogap} and Theorem~\ref{thm:gapone} simpler?
\end{question}


 



Furthermore, there is a notion of \newword{order-chain} polytopes of a poset $P$, which generalizes the order polytope and the chain polytope and creates a class of polytopes that sort of interpolates between the two \cite{Hibi2015OrderchainP}, \cite{FANG2016267}. It is natural to ask if our notion of crossing-number can also provide a bound on the difference between the number of facets of such polytopes:

\begin{question}
 Can the crossing number be used to come up with the facet number gaps between various order-chain polytopes?   
\end{question}

\subsection*{Acknowledgments}
This research was carried out primarily during the $7331$ Combinatorics course, held during the $2024$ spring semester at Texas State University. The authors appreciate the support from Texas State University for providing support and a great working environment. The authors also thank Takayuki Hibi for introducing the authors to the topic.


\bibliography{poset_bib}

\begin{thebibliography}{1}

\bibitem{ahmad2023order}
{\sc I.~Ahmad, G.~Fourier, and M.~Joswig}, {\em Order and chain polytopes of maximal ranked posets}, arXiv preprint arXiv:2309.01626,  (2023).

\bibitem{FANG2016267}
{\sc X.~Fang and G.~Fourier}, {\em Marked chain-order polytopes}, European Journal of Combinatorics, 58 (2016), pp.~267--282.

\bibitem{freijhollanti2024fvector}
{\sc R.~Freij-Hollanti and T.~Lundström}, {\em $f$-vector inequalities for order and chain polytopes}, arXiv preprint arXiv:2312.13890,  (2024).

\bibitem{HIBILI16}
{\sc T.~HIBI and N.~LI}, {\em Unimodular equivalence of order and chain polytopes}, Mathematica Scandinavica, 118 (2016), pp.~5--12.

\bibitem{Hibi2015OrderchainP}
{\sc T.~Hibi, N.~Li, T.~X. Li, L.~Mu, and A.~Tsuchiya}, {\em Order-chain polytopes}, Ars Math. Contemp., 16 (2015), pp.~299--317.

\bibitem{HIBI17}
{\sc T.~Hibi, N.~Li, Y.~Sahara, and A.~Shikama}, {\em The numbers of edges of the order polytope and the chain polytope of a finite partially ordered set}, Discrete Mathematics, 340 (2017), pp.~991--994.

\bibitem{stanley1986two}
{\sc R.~P. Stanley}, {\em Two poset polytopes}, Discrete \& Computational Geometry, 1 (1986), pp.~9--23.

\end{thebibliography}
\bibliographystyle{siam}

\end{document}